\documentclass[12pt,reqno] {amsart}
\usepackage[mathscr]{euscript}
\usepackage{amsthm}
\usepackage{amssymb}
\usepackage{amsmath}%Gamma
\usepackage{verbatim}
\usepackage{enumerate}
\usepackage{amscd}
\usepackage{  latexsym }
\numberwithin{equation}{section}
\newcommand{\cala}{\mathcal A} 
\newcommand{\calb}{\mathcal B} 
\newcommand{\calu}{\mathcal U} 
\newcommand{\osn}{\mathcal{ O}_s(N)} 
\newcommand{\tsn}{\mathcal{ T}_s(N)} 
\newcommand{\sumti}{\sum_{\langle t_i \rangle}} 
\newcommand{\sumpi}{\sum_{\langle p_i \rangle}}
\newcommand{\phiquotient}{\left(\prod_{1 \le i \le s} \phi(t_i)\right)/\phi(t)}
\swapnumbers
\theoremstyle{plain}
 \newtheorem{theorem}{Theorem}[section]
\newtheorem{main}[theorem]{Main Theorem}
\newtheorem{lemma}[theorem]{Lemma}
\newtheorem{corollary}[theorem]{Corollary} 
\newtheorem{proposition}[theorem]{Proposition}
\theoremstyle{remark}
\newtheorem{remark}[theorem]{Remark}

\begin{document} 
\date{\today}
\title{Many Frobenius complements \\ have even order}  

\author{Ron Brown}
\address{Department of Mathematics\\University of Hawaii\\2565 McCarthy Mall\\
Honolulu, HI 96822}
\email{ron@math.hawaii.edu}
 
 \keywords{ Frobenius group, Frobenius complement, Sylow-cyclic group, Frobenius triple, distribution of primes}
\subjclass[2010]{Primary:  20D60; Secondary:  11N13, 20E99}  
 
\begin{abstract}
The theory of Frobenius groups with Frobenius complements of even order largely reduces to tractable algebraic number theory.  If we consider only Frobenius complements with an upper bound $s$ on the number of distinct primes dividing the order of their commutator subgroups, then the proportion of these with odd order is less than $1/2^{s}$. A positive lower bound   is also given.
   
\end{abstract} 
\maketitle
%\newpage
%\tableofcontents

%\pagestyle{empty}
%\renewcommand{\baselinestretch}{1.5} 

%\thispagestyle{empty}

\section{Introduction and the Main Theorem}\label{s: intro}

Frobenius groups play a significant role in the theory of finite groups, and, in particular, in the study of simple groups. For example it is not unusual in finite group theory that a proof analyzing a minimal counterexample will end up considering Frobenius groups.   Every Frobenius group $G$ is a semidirect product $K \rtimes H$ where $K$ is a canonical normal subgroup of $G$ (the \emph{Frobenius kernel} of $G$) and $H$ is a subgroup of $G$ (a \emph{Frobenius complement} of $G$) \cite[35.25(1)]{A}.  The theory of Frobenius groups with abelian Frobenius kernel largely reduces to algebraic number theory and indeed to a tractable part of algebraic number theory: the study of unramified primes in abelian extensions of the field of rational numbers  \cite{Br}.  For example it is fairly easy to count the exact number of isomorphism classes of such groups of order less than $10^{6}$; there are $569,342$ of them \cite[p. 85]{Br}.  It is natural then to ask if it is common for Frobenius groups to have abelian kernel.

QUESTION.  Do almost all (or even a positive proportion of) Frobenius groups have abelian Frobenius kernel?

When the complement of a Frobenius group has even order, then it is known that the kernel is abelian \cite[Theorem 3.4A]{DM}.  About $88.8\%$ of Frobenius complements of order at most $10^6$ have even order \cite[p. 54]{Br}.   

CONJECTURE. Almost all (isomorphism classes of)    Frobenius complements have even order.

 In fact,  all Frobenius complements have even order except for some of those which are $\mathbb{Z}$-groups, i.e., groups all of whose Sylow subgroups are cyclic \cite[Theorem 1.4]{Br}.  We follow Lam in calling these groups
\textit{Sylow-cyclic} groups \cite{L}.  Our Main Theorem below would appear to support the stronger conjecture that almost all Sylow-cyclic Frobenius complements have even order.

By the \textit{breadth} of a Sylow-cyclic Frobenius complement we mean the number of distinct primes dividing the order of its commutator subgroup.  In this paper we will show that most of the members of the family of Sylow-cyclic Frobenius complements whose breadth is bounded above by some large number have even order.  More precisely we prove the 

\begin{main}\label{main}  
Suppose that $s \ge 1$.  Then for all sufficiently large $N$, the proportion of all Sylow-cyclic Frobenius complements of order at most $N$ and breadth at most $s$ which have odd order is less than $1/2^{s }$. 
\end{main}

The required size of $N$ will of course depend on the choice of $s$.  The above ``proportion" really refers to the proportion of isomorphism classes of 
Sylow-cyclic Frobenius complements.  We will often omit the phrase ``isomorphism classes of" below, just as we did in the statement of the Main Theorem.  

%The reader will find it easy to sharpen the above result (as well as some of the preliminary results).  However it seemed to us that in the absence of significant %applications, the complications involved in obtaining these improvements were not justified at this time.

 In the next section we establish some formulas for the number $\tsn$ of Sylow-cyclic Frobenius complements of breadth $s$ and order at most $N$ and the number $\osn$ of these of odd order.  After establishing some preliminary lemmas in Section 3 we give a lower bound for $\tsn$  (Section 4) and an upper bound for $\osn$ (Section 5).
These are used in Section \ref{s: mainproof} to   show that
\[
\rho_s(N) :=\sum_{r \le s} \mathcal{O}_r(N) \Big/\sum_{r \le s} \mathcal{T}_r(N)   
\]
 is less than $1/2^s$ for sufficiently large $N$. In Section \ref{s: morebounds} we sketch a computation of the exact proportion of Sylow-cyclic Frobenius complements with breadth at most $s$ which have odd order and use this result to provide a   lower bound and an improved upper bound for this proportion.

Some conventions with our notation will prove convenient.  Lower case Roman letters (with the exception of $r$ and $x$), with or without subscripts, will always denote positive integers; $r$ will always denote a nonnegative integer.  The lower case $p$, with or without subscripts, will always denote a prime.  The notation $t_0$ wll denote the product of the distinct primes dividing $t$ (so, for example, $12_0=6$ and $1_0=1)$.  The least common multiple of $a_1,\cdots,a_s$ will be denoted by $[a_1,\cdots,a_s]$.  $|\calb|$ denotes the number of elements in the finite set $\calb$. $\phi$ of course denotes Euler's $\phi$-function.

\textit{Throughout the paper $s$ will denote a fixed positive integer.}  The letters $A, B, C,  D, E, F, F', G, H, H', I, I', J, J'$ will denote various constants, some depending on the choice of $s$.  Finally, $N$ will always be assumed to be a positive integer large enough that $\log\log\log N > 2$ and $N> (\log\log N)^{6s + 6}$.

The prerequisites from \cite{Br} for reading this paper are modest; specifically we use 5.1 (A), (B) and (C) (definitions), 5.2 (A), 5.3 (A) (in the statement of which the numbers $2$ and $3$ were unfortunately transposed), and 11.1 (A) and (B1); all of this is largely self-contained. Sylow-cyclic Frobenius complements are also called ``$1$-complements" in \cite{Br}.

\section{Two counting formulas}\label{s: counting}

The isomorphism classes of Sylow-cyclic Frobenius complements of order at most $N$ and breadth $s$ correspond bijectively   to the ``proper Frobenius triples of order at most $N$ and breadth $s$", i.e.,
to triples $(m,n, \langle r + m\mathbb{Z} \rangle)$ (where $\langle r + m\mathbb{Z} \rangle$ is the multiplicative group generated by a unit $r+m\mathbb{Z}$ of the ring 
$\mathbb{Z}/m\mathbb{Z}$) such that $n>1$, $m$ is relatively prime to   $n$ and to  $r-1$ (so $2 \nmid m$),  $r^{n/n_0} \equiv 1 \pmod{m} $, $mn \le N$ and  $m$ is divisible by exactly $s$ distinct primes 
 \cite[Theorem 5.2 A and Lemma 5.3 A]{Br}. %Note that for  such triples the multiplicative order of $r + m \mathbb{Z}$ divides $n/n_0$.    
Each such triple uniquely defines an ordered pair
\begin{equation}\label{opair}
  (k,\{ (p_1^{a_1}, t_1), \cdots, (p_s^{a_s},t_s) \} )
\end{equation}
where $p_1^{a_1} \cdots p_s^{a_s}$  is the prime factorization of  $m$, each $t_i$ is the multiplicative order of $r+ p_i\mathbb{Z}$ in the ring 
$\mathbb{Z}/p_i\mathbb{Z}$,  and $k = n/t t_0$ where $t=[t_1,\cdots,t_s]$; by \cite[Proposition 11.1 (A)]{Br} these ordered pairs satisfy for all $i \le s$ and $j \le s$:

(1) $t_i > 1 \mbox{ and } p_i \equiv 1 \pmod{t_i}$;

(2) $p_j \nmid t_i \mbox{ and if } i \ne j, \mbox{ then } p_i \ne p_j$;

(3) $p_1^{a_1} \cdots p_s^{a_s}  [t_1,\cdots,t_s][t_1,\cdots,t_s]_0 \le N$; and

(4) $p_j \nmid k$ and $p_1^{a_1} \cdots p_s^{a_s} [t_1,\cdots,t_s] [t_1,\cdots,t_s]_0 k \le N$.    \newline 
 Note that $k$ above is an integer since  $r^{n/n_0} \equiv 1 \pmod{m} $.

On the other hand each ordered pair of the form (\ref{opair}) satisfying the above four conditions  arises as above from exactly
 $\frac{1}{\phi([t_1, \cdots, t_s])} \prod_{1 \le i \le s}\phi (t_i)$ proper Frobenius triples of order at most $N$ and breadth $s$ \cite[Proposition 11.1  (B1)]{Br}.  Let $\mathcal{D}$ denote the set of such ordered pairs. 

Each   ordered pair (\ref{opair}) in $\mathcal{D}$ gives rise to exactly
$s!$ ordered pairs
\begin{equation}\label{s+tuple}
(k,  ((p_1^{a_1}, t_1), \cdots, (p_s^{a_s},t_s))  )
\end{equation}
satisfying the four conditions above.   For any $t>1$ let $\mathcal{D}_t$ denote the set of ordered pairs (\ref{s+tuple}) with 
$t=[t_1, \cdots, t_s]$ and also satisfying the four conditions above.

For any $t>1$ let  $\calu_N(t)$ denote the set of $s$-tuples
\begin{equation}\label{2.1}
\tau = ((p_1^{a_1},t_1), \cdots, (p_s^{a_s},t_s))
\end{equation}
with $t= [t_1, \cdots, t_s]$  and also satisfying  conditions (1), (2), and (3) above; for any such $\tau \in \calu_N(t)$ set 
\[
\phi(\tau) = \frac{1}{\phi(t)} \prod_{1 \le i \le s}\phi (t_i)
\]
and 
\[
K(\tau) = \{ k:  (k,\tau) \in \mathcal{D}_t \}.
\]

If $\tau^*$ is the set of all $s$ coordinates of $\tau \in \calu_N(t)$, then we also set $\phi(\tau^*)=\phi(\tau)$.
 
Note that if $\calu_N(t)$ is nonempty, so that it contains some $\tau$ as in (\ref{2.1}), then
\[
t^2 \le t t_1 \cdots t_s 
\le tt_0p_1^{a_1} \cdots p_s^{a_s} \le N,
\]
so $t \le \sqrt{N}$.

Combining these several observations we can give formulas for $\tsn$ and $\osn$.

\begin{proposition}\label{prop1}
 
\[
s!\tsn 
 = \sum_{2 \le t \le \sqrt{N}}  \quad  \sum_{\tau \in \calu_N(t)}   | K(\tau)| \phi(\tau).
\]
\end{proposition}

\begin{proof} 
 \[ 
s!\tsn = s! \sum_{(k,\tau)\in \mathcal{D}} \phi(\tau) = \sum_{2 \le t \le \sqrt{N}} \quad \sum_{(k,\tau) \in \mathcal{D}_t} \phi(\tau)
\]
\[
 = \sum_{2 \le t \le \sqrt{N}} \quad \sum_{\tau \in \calu_N(t)} |K(\tau)|\phi(\tau).
\]
\end{proof}

A similar argument yields the next
   
\begin{proposition}\label{prop2}
 
\[
s!\osn = \sum_ {\substack{ 2 \le t \le \sqrt{N} \\ 2 \nmid t} } \quad \sum_{\tau \in \calu_N(t)}   |\{ k\in K(\tau): 2 \nmid k \}|  \phi(\tau).
\]
\end{proposition}

\section{Some preliminary lemmas}\label{s: prelims}

We continue to fix $s \ge 1$.  For any $t> 1$ we set 
\[ 
\Gamma(t) = \{ (a_1, \cdots, a_s): [a_1, \cdots,  a_s] = t\}
\] 
and
\[
\Gamma^*(t) = \{ (a_1, \cdots, a_s)\in \Gamma(t): a_i> 1  \mbox{ for all } i \le s\}.
\]
We also set $\gamma(t) = |\Gamma(t)|$,   $\gamma^*(t) = |\Gamma^*(t)|$, and $\gamma(1) = 1$.  

\begin{lemma}\label{Lem3.1}
 For all $t>1$ and $r, p, m,$ and  $n$ we have:

(A) if $m$ and $n$ are relatively prime, then $\gamma(mn) = \gamma(m)\gamma(n)$;

(B) if $2\nmid t$, then $\gamma^*(2^rt) \ge \gamma (2^r) \gamma^*(t)$;  

(C) $\gamma(p^r) = (r+1)^s -r^s$; and

(D) $\gamma^*(p^m) = m^s -(m-1)^s$.

\end{lemma}

\begin{proof}
 Parts (B) and (C) follow trivially from the definition of $\gamma(1)$ if $r=0$ so suppose $r>0$.  Similarly we may suppose that $m>1$ and $n>1$.
We have inverse maps $\alpha:\Gamma(m) \times \Gamma(n) \longrightarrow \Gamma(mn)  $ and $\beta:\Gamma(mn) \longrightarrow \Gamma(m) \times    \Gamma(n)$ with
\[
 \alpha((a_1,\cdots,a_s),(b_1,\cdots,b_s)) = (a_1b_1,\cdots,a_sb_s) 
\] 
and
\[
\beta(a_1,\cdots,a_s) = (((a_1,m), \cdots,  (a_s, m)), ((a_1,n),  \cdots, (a_s, n)))
\]
where in the above display we have let ``$(a,b)$" denote the greatest common divisor of $a$ and $b$. This implies part (A) above; part (B) follows since we can restrict the map $\alpha$ to an injective map $\Gamma(2^r) \times \Gamma^*(t) \longrightarrow \Gamma^*(2^rt)$ when $m=2^r$ and $n = t$.

  Next, an element of $\Gamma(p^r)$ will have the entry $p^r$ in some coordinate; for any $1 \le j \le s$ there are $\binom{s}{j}$ ways of filling $j$ coordinates with the power $p^r$; the remaining $s-j$ coordinates can be filled with any combination of the powers $p^{k-1}$ where $1 \le k \le r$.  Therefore the number of elements of $\Gamma(p^r)$
is 
\[
 \sum_{1 \le j \le s} \binom{s}{j} r^{s-j} 1^j = (r+1)^s-r^s.
\]

Similarly, an element of $\Gamma^*(p^m) $ must have entry $p^m$ in $j$ coordinates where $1 \le j \le s$ and the remaining $s-j$ coordinates can be filled with any
$p^i$ where $1  \le i <m$, so $\gamma^*(p^m)= \sum_{1 \le j \le s} \binom{s}{j} (m-1)^{s-j} 1^j = m^s - (m-1)^s$.
\end{proof}.

\begin{remark}
Assertions (A) and (C) of the preceding lemma give a way of computing $\gamma(t)$ from the prime factorization of $t$.  Let us in this remark write 
$\gamma_s$ instead of $\gamma$ to indicate the dependence of the function $\gamma$ on the choice of $s$.  Then we can also compute the values of $\gamma^*(t)$ using the formula 
\[
\gamma^* = \sum_{r=0}^{s-1} (-1)^r \binom{s}{r} \gamma_{s-r} .
\]

\end{remark}

We will set $f(t):= 1/(t t_0 \phi(t))$.  Note that $f(t) = 1/(t^2\phi(t_0))$.
\begin{lemma}\label{D}
The series
\[
  \sum_{2 \le t < \infty} f(t)\gamma^*(t) (\log t)^s
\]

converges. 
\end{lemma}
\begin{proof}
There is a constant $A>0$ (depending only on $s$) such that the number $d(m)$ of positive divisors of any positive integer $m$ is at most $Am^{1/(10s)}$ \cite[Theorem 315]{HW}.  Then 
\[
 \sum_{t = 2}^{\infty} f(t)\gamma^*(t) (\log t)^s \le \sum_{t = 2}^{\infty} f(t)d(t)^s (\log t)^s \le \sum_{t = 2}^{\infty} \frac{A^s t^{1/10} (\log t)^s}{t^2 \phi(t_0)} < \infty.
\]
\end{proof} 
 \begin{corollary}\label{corD}
The series $\sum_{ \substack { 2 \le t < \infty\\2\nmid t}} f(t)\gamma^*(t)$ and $\sum_{  2 \le t < \infty} f(t)\gamma^*(t)$ both
converge.
\end{corollary} 

 In   Theorem \ref{exact} below we   give the exact proportion of Sylow-cyclic Frobenius complements of breadth at most $s$ which have odd order in terms of the convergent series of the above corollary.

\begin{lemma}\label{triplelog}
For integers $t 
\ge 2$ and $u \ge 2^3$ we have 
\[
\phi(t) \sum_{ \substack{ p \equiv 1 \pmod{t}\\ p\mid u}} \frac{1}{p} \le 3 + 2 \log\log(\log u / \log 2).
\]
\end{lemma}

\begin{proof}
Let $q_1^{c_1} \cdots q_r^{c_r}$ be the prime factorization of $u$, so that $u \ge 2^r$ and hence $r \le  \log u / \log 2$.  For each $j \ge 1$ let $p_j$ denote 
the $j$-th prime congruent to $1$ modulo $t$, so $p_j \ge 1 + jt$ and so $j  \le p_j/t$.  By the  Brun-Titchmarsh theorem \cite{BT} if $j>1$, then  the number $j$ of primes congruent to $1$ modulo $t$ and at most equal to $p_j$ satisfies
\[
j \le \frac{2p_j}{\phi(t)\log(p_j/t)} \le \frac{2p_j}{\phi(t)\log j},
\]
so 
\begin{equation}\label{3.1}
 p_j  \ge \frac{1}{2} \phi(t) j \log j,
\end{equation}
 which is of course also true if $j=1$. The lemma is trivally true if $r < 4$ since by our choice of $u$ 
we have $\log\log(\log u / \log 2) \ge 0$.  Hence suppose that $r \ge 4$.  Then the inequality (\ref{3.1}) implies that
\[
\phi(t) \sum_{ \substack{ p \equiv 1 \pmod{t}\\ p\mid u}} \frac{1}{p} \le \phi(t) \left( \frac{1}{p_1}+\frac{1}{p_2}+\frac{1}{p_3}\right) + 2\int_3^r \frac{dx}{x\log x}
\]
\[
 \le 3 + 2\log\log r \le 3 + 2 \log\log(\log u /\log 2).
\]
\end{proof}

Note.  If $u \le \sqrt{N}$, then  by the above theorem
\[
 \phi(t) \big(\sum_{ \substack{ p \equiv 1 \pmod{t}\\ p\mid u}} \frac{1}{p}\big) \le 3 + 2 \log\log(\log \sqrt{N}/ \log 2) \le 3 + 2\log\log\log N .
\]

The next lemma is simply a variant of a special case of a lemma of K. K. Norton.

\begin{lemma}\label{carl}
Suppose $t>1$.  Then there is a constant $C>2$ (independent of $N$ and $t$) with 
\[ 
\left|\phi(t) \sum_{ \substack{ p \equiv 1 \pmod{t}\\ p \le N}} \frac{1}{p} - \log\log N \right| \le C \log t.
\]
\end{lemma}

\begin{proof}
By \cite[Lemma 6.3]{N} there is a constant $B$ (independent of $N$ and $t$) such that 
\[
\left|\sum_{ \substack{ p \equiv 1 \pmod{t}\\ p \le N}} \frac{1}{p} -
\frac{\log\log N}{\phi(t)}   \right| \le B \frac{\log (3t)}{\phi(t)}
\] 
(in Norton's lemma set $L=\{ 1 \}, k = t$, and $x=N$).  The lemma follows immediately since $t \ge 2$. 
\end{proof}

\section{Lower bounds for $\tsn$}\label{s: lower}

It will be convenient to use in this section some abbreviations involving summation notation.  In the expression
\[
\boldsymbol{\sum} = \sum_t\sumti\sumpi\sum_{p_s}
\]
we let $\sum_t$ sum over all $t$ with $1< t \le \log\log N$; $\sumti$ will abbreviate $\sum_{t_1} \cdots \sum_{t_s}$ where for each $i<s$, $\sum_{t_i}$ sums over all $t_i>1$ with $t_i \mid t$ and $\sum_{t_s}$ sums over all $t_s > 1$ with         $[t_1, \cdots, t_s] = t$.  (Thus if $s=1$ we have $\sum_t = \sum_t \sumti$.)
Also $\sumpi$ abbreviates $\sum_{p_1} \cdots \sum_{p_{s-1}}$ where each $\sum_{p_i}$ sums over all primes $p_i \equiv 1 \pmod{t_i}$ with $p_i \nmid t$, $p_i \nmid  p_1 \cdots p_{i-1}$, and $p_i \le N^{1/(3s)}$; and, finally, where $\sum_{p_s}$ sums
over all $p_s\equiv 1 \pmod{t_s}$ with $p_s \nmid t$, $p_s \nmid  p_1 \cdots p_{s-1}$ and $t t_0 p_1 \cdots p_s \le \sqrt{N}$.  We are interpreting empty products (such as $ p_1 \cdots p_{s-1}$ if $s=1$) as equal to $1$.  Thus, for example, if $s=1$ then  $\sum_{p_s}$ sums
over all $p_s\equiv 1 \pmod{t_s}$  (so consequently $p_s \nmid t$)  and $t t_0   p_s \le \sqrt{N}$, and $\boldsymbol{\sum} = \sum_t \sum_{p_s}$. Note that for all $s$, 
\[
\gamma^*(t) =  \sumti 1   \mbox{ for all }  t>1.
\]

It must be kept in mind that this notation depends heavily on context.  Thus, for example, the precise meaning of $\sum_{t_s}$ depends on the previous choices of the  parameters $t, t_1, \cdots, t_{s-1}$,  and similarly 
$\sum_{p_i}$ depends on the previous choices of $t, t_i, p_1, \cdots, p_{i-1}$ (as well as $N$).  For some choices of $t, t_1, \cdots, t_s, p_1, \cdots, p_{s-1}$ it could be the case that $\sum_{p_s}$ is an empty sum.  If not, and $t, t_1, \cdots, t_s, p_1, \cdots, p_s$ is a set of parameters for a term of the summation 
$\boldsymbol{\sum}$, then $1<t 
\le \log\log N < \sqrt{N}$ and $((p_1,t_1), \cdots, (p_s,t_s)) \in \calu_N(t)$.  It therefore follows from
Proposition \ref{prop1} that
\begin{equation}\label{4.1}
s!\tsn \ge \boldsymbol{\sum} \frac{\prod_{1 \le i \le s} \phi(t_i)}{\phi(t)} \left| K(((p_1,t_1), \cdots, (p_s,t_s)))\right|.
\end{equation}

\begin{lemma}
Suppose that $t$, $N$ and the $s$-tuple $\tau=((p_1,t_1), \cdots, (p_s,t_s))$ appear in the inequality (\ref{4.1}) (with $t=[t_1, \cdots, t_s]$).  Set $\Delta = N/(t t_0 p_1 \cdots p_s)$.  Then
$|K(\tau)| \ge \Delta \prod_{1 \le i \le s} (1-\frac{1}{p_i}) \quad - \quad 2^{s }$.
\end{lemma}

\begin{proof}
For any subset $\alpha$ of $\{ 1, 2, \cdots,s\}$ let $P_{\alpha} =\prod_{i \in \alpha} p_i$ (so $P_{\emptyset} = 1$).  For any $r\le s$  let $\cala_r$ denote the set of subsets of $\{ 1, \cdots,s\}$ of order $r$.
An inclusion/exclusion argument \cite[Theorem 261]{HW} shows that 
\[
|K(\tau)| = \sum_{r=0}^s (-1)^r \sum_{\alpha \in \cala_r} \left[ \frac{\Delta}{P_{\alpha}}\right]
\]
\[
 = \sum_{r=0}^s (-1)^r \sum_{\alpha \in \cala_r}   \frac{\Delta}{P_{\alpha}}\quad -\quad \sum_{r=0}^s (-1)^r  \sum_{\alpha \in \cala_r} \big( \frac{\Delta}{P_{\alpha}} - \left[ \frac{\Delta}{P_{\alpha}}\right]\big)
\]
\[
\ge \Delta\prod_{i=1}^s(1-\frac{1}{p_i}) \quad -\quad \sum_{r=0}^s    \sum_{\alpha \in \cala_r} 1 \ge \Delta\prod_{i=1}^s(1-\frac{1}{p_i})  - 2^{s }.
\]
\end{proof}.

\begin{lemma}
Let $i \le s$. 
For all choices of $t, t_1, \cdots , t_i, p_1, \cdots,p_{i-1}$ (or just of $t$ and $t_1$ if $i=1$) in display  
 (\ref{4.1}) we have 
\begin{equation}\label{4.2}
\sum_{p_i} \phi(t_i)\left(\frac{1}{p_i}-\frac{1}{p_i^2}\right) \quad\ge\quad \log\log N - E \log \log \log N
\end{equation}
for a positive constant $E$ (depending only on $s$). 
\end{lemma} 
\begin{proof}
By our choice of $N$ in \S\ref{s: intro} we have $\sqrt{N}/(tt_0p_1\cdots p_{s-1}) \ge N^{1/(3s)}$.  Thus  
\begin{equation}\label{4.3}
%\begin{split}
\sum_{p_i} \phi(t_i)\left(\frac{1}{p_i}-\frac{1}{p_i^2}\right)  
%\\ 
\ge \sum_{p_i \le N^{1/(3s)}}  \frac{\phi(t_i)}{p_i} - \sum_{p_i | tp_1 \cdots p_{s-1}}\frac{\phi(t_i)}{p_i} -  \sum_{p_i \le N} \frac{\phi(t_i)}{p_i^2}
%\end{split}
\end{equation}

where each sum above is only over primes $p_i \equiv 1  \pmod{t_i}$.

First, 
\[
 \sum_{\substack{ p_i \le N\\ p_i \equiv 1 \pmod{t_i}}} \frac{\phi(t_i)}{p_i^2} \le \sum_{k=1}^{\infty} \frac{t_i}{(1+kt_i)^2} \le \sum_{k=1}^{\infty} \frac{1}{k^2} 
= \frac{\pi^2}{6}.
\] 

Next, if the number of prime divisors of $tp_1 \cdots p_{s-1}$ is less than $4$. then
\begin{equation}\label{4.4} 
\sum_{\substack{p_i | tp_1 \cdots p_{s-1}\\p_i \equiv 1 \pmod{t_i}}}\frac{\phi(t_i)} {p_i} \le 3 + 2\log\log\log N.
\end{equation}
If the number is greater than 3, then
\[
2^3 \le tp_1 \cdots p_{s-1} \le (\log\log N)(N^{1/(3s)})^{s-1} \le \sqrt{N}
\]
so by Lemma \ref{triplelog} (and the Note following it), the   inequality (\ref{4.4})
 still holds.

Finally, by Lemma \ref{carl} and the choice of $t$ we have
\[
\sum_{\substack{ p_i \le N^{1/(3s)}\\ p_i \equiv 1 \pmod{t_i}}} \frac{\phi(t_i)}{p_i} \ge \log\log N^{1/(3s)} - C \log t_i
\]
\[
\ge \log\log N - \log 3s -C\log t \ge \log\log N - D \log\log\log N
\]
for a constant $D$ depending only on $s$.

Applying the inequalities of the last three paragraphs to the inequality (\ref{4.3}) we deduce that
\[
\sum_{p_i} \phi(t_i)\left(\frac{1}{p_i}-\frac{1}{p_i^2}\right) \ge \log\log N - D \log \log \log N - \frac{\pi^2}{6} -3 -2\log\log\log N
\]
\[
\ge \log\log N - E \log\log\log N
\]
where the positive constant $E$ depends only on s.
\end{proof}

Observe that since $N \ge (\log\log N)^{6s+6}$, therefore
\[
\boldsymbol{\sum} \phiquotient  \le \sum_t \sumti\sumpi \frac{\sqrt{N}}{ t t_0 p_1 \cdots p_{s-1}}  \phiquotient
\]
\[
\le \sqrt{N}\sum_t \sumti\sumpi 1 \le \sqrt{N}(\log\log N)^{s+1}(N^{1/(3s)})^{s-1}  \le N.
\]
 
We can now give our lower bound for $\tsn$.

\begin{theorem}\label{lower}
There is a constant $I$ depending only on $s$ such that $s!\tsn $ is greater than or equal to 
\[
  N\left(\sum_t f(t)\gamma^*(t)\right)(\log \log N)^s \quad   - \quad    I
 N(\log\log N)^{s-1}\log\log\log N.
\]
\end{theorem}

\begin{proof}
 Applying the last paragraph and the previous two lemmas to the terms of display (\ref{4.1}) we have 
\[
s!\tsn \ge \boldsymbol{\sum}\frac{1}{\phi(t)} \left(\prod_{1 \le i \le s}\phi(t_i)\right)\left(- 2^{s} +  \frac{N}{t t_0 p_1\cdots p_s} \prod_{1 \le i \le s} \left(1-\frac{1}{p_i}\right)   \right) 
\]
\[
\ge - 2^{s}N + N \sum_t f(t)\sumti\sumpi\sum_{p_s}\prod_{1 \le i \le s}(\phi(t_i)(\frac{1}{p_i} - \frac{1}{p_i^2})) 
\]
\[
\ge - 2^{s}N + N\left(\sum_t f(t)\gamma^*(t)\right)(\log\log N   -   E  \log\log\log N)^s  
\]
\[
\ge N\left(\sum_t f(t)\gamma^*(t)\right)(\log\log N )^s \quad -\quad  I N (\log\log N)^{s-1}\log\log\log N
\]
for a constant $I>0$ which depends only on $s$.  (Recall that $\sum_t f(t)\gamma^*(t)$ converges by  Corollary \ref{corD}.)
\end{proof}

\section{Upper bounds for $\osn$}\label{upper}

We introduce abbreviations for summation notation to be used only in this section.  We will consider expressions 
\[
\boldsymbol{\sum} = \sum_t\sumti\sum_{\langle p_i,a_i\rangle} \sum_{p_s}\sum_{a_s}
\]
where $\sum_t$ sums over all odd $t$ with $3 \le t \le \sqrt{N}$; $\sumti$ has the same meaning as in \S \ref{s: lower}  (so again $\sumti 1 = \gamma^*(t)$); 
$\sum_{\langle p_i,a_i\rangle}$ abbreviates $\sum_{p_1}\sum_{a_1} \cdots \sum_{p_{s-1}}\sum_{a_{s-1}}$ where for each $i\le s$, $\sum_{p_i}$ is the sum over all $p_i \equiv 1 \pmod{t_i}$ with $p_i \le N$ and $p_i \nmid t p_1 \cdots p_{i-1}$; and for $i<s$, $\sum_{a_i}$ is the sum over all $a_i$ with $p_i^{a_i} \le N$, while
$\sum_{a_s}$ is the sum over all $a_s$ with $t t_0 p_1^{a_1} \cdots p_s^{a_s} \le N$.

Then arguing as in the proof of the inequality (\ref{4.1}) we can see that Proposition \ref{prop2} implies that
\begin{equation}\label{5.1}
s!\osn \\ \le   \boldsymbol{\sum} \left( \frac{\prod_{1\le i \le s} \phi(t_i)}{\phi(t)}\right) \left| \{ k \in K(((p_1^{a_1}, t_1), \cdots , (p_s^{a_s},t_s))): 2 \nmid k\}\right|.
\end{equation}
After all, any sequence $t_1, \cdots , t_s, p_1, a_1, \cdots, p_s, a_s$ which satisfies   conditions (1), (2) and (3) in Section \ref{s: counting} and has $2\nmid t_i$ for all $i \le s$ will, with $t=[t_1, \cdots , t_s]$, be a sequence of parameters in the summation $\boldsymbol{\sum}$.   

In this section we set $p=2$ so $p_0  =2$.  Consider an $s$-tuple $\tau = ((p_1^{a_1},t_1), \cdots, (p_s^{a_s},t_s))$ appearing in the inequality (\ref{5.1}). For any $r\le s+1$ let $\calb_r$ denote the set of subsets of $\{ 0,1,\cdots , s\}$ of order $r$; for any subset $\beta$ of $\{0,1,\cdots,s\}$ let $P_{\beta} = \prod_{r \in \beta} p_r$.  Also set $\Delta = N/(tt_0p_1^{a_1} \cdots p_s^{a_s})$.  Arguing much as in \S  \ref{s: lower}  we have  
\[
\left| \{ k  \in K(\tau): 2 \nmid k\}\right|= \sum_{0 \le r \le s+1}\sum_{\beta \in \calb_r} (-1)^r\left[\frac{\Delta}{P_{\beta}}\right]
\]
\[
\le \sum_{0 \le r \le s+1} (-1)^r\sum_{\beta \in \calb_r} \frac{\Delta}{P_{\beta}}  +        \sum_{ \substack{ 0 \le r  \le s+1\\2 \nmid r}}  \:  \sum_{\beta \in \calb_r}\left( \frac{\Delta}{P_{\beta}}-         
 \left[\frac{\Delta}{P_{\beta}}\right]\right)
\]
\[
\le \Delta \prod_{0 \le r \le s } \left( 1 - \frac{1}{p_r}\right) + \sum_{ \substack{ 0 \le r  \le s+1\\2 \nmid r}} | \calb_r| = \frac{\Delta}{2} \prod_{1 \le r \le s} \left( 1 - \frac{1}{p_r}\right) + 2^s
\]
since $2^s$ is the number of subsets of $\{ 0,1,\cdots s\}$ of odd order.  Hence

\begin{equation}\label{5.2}
s!\osn \le \frac{N}{2} \boldsymbol{\sum} f(t)   \prod_{1 \le i \le s}\phi(t_i)\left(  \frac{1}{p_i^{a_i}}(1 - \frac{1}{p_i }) \right)  \quad +\quad 2^s  \boldsymbol{\sum}\frac{\prod_{1\le i \le s} \phi(t_i)}{\phi(t)}.
\end{equation}

\begin{lemma}
For some constant $G$ depending only on $s$ we have 
\[
\boldsymbol{\sum} \phiquotient \le GN(\log\log N)^{s-1}.
\]
\end{lemma}

\begin{proof}
Given $t,t_1, \cdots, t_s,p_1,a_1,\cdots, p_{s },a_{s }$ as in the summation notation $\boldsymbol{\sum}$ for this section we have  $p_s^{a_s} \equiv 1 \pmod{t_s}$, so
\[
\sum_{p_s}\sum_{a_s}\phi(t_s) \le  \phi(t_s) \left|\{k: 1+kt_s \le N/( t t_0 p_1^{a_1}\cdots p_{s-1}^{a_{s-1}})\}\right|
\]
\[ 
\le (\phi(t_s)/t_s)  N/(t t_0 p_1^{a_1}\cdots p_{s-1}^{a_{s-1}}) \le N/(t t_0 p_1^{a_1}\cdots p_{s-1}^{a_{s-1}}).
\]
 (Here, if $s=1$, we regard the expression $ p_1^{a_1}\cdots p_{s-1}^{a_{s-1}}$ as an empty product, so it is equal to $1$.)

Also for all $i<s$, 
\[
\phi(t_i)\sum_{p_i} \sum_{a_i} \frac{1}{p_i^{a_i}} \le \phi(t_i)\sum_{p_i}\sum_{k=1}^{\infty} \frac{1}{p_i^k}
\]
\[
 \le \phi(t_i)\sum_{p_i}\frac{2}{p_i} \le 2 \phi(t_i) \sum_{\substack{ p_i \le N \\ p_i \equiv 1 \pmod{t_i}}} \frac{1}{p_i} \le 2 \log\log N + 2C\log t_i
\]
by Lemma \ref{carl}.  Therefore since $2\log t > 1$ for $t>1$, 
\[
\boldsymbol{\sum} \frac{\prod_{1\le i \le s} \phi(t_i)}{\phi(t)}
\]
\[
 \le \sum_t \frac{1}{\phi(t)} \sumti \left( \sum_{p_1}\sum_{a_1}\frac{\phi(t_1)}{p_1} \cdots \sum_{p_{s-1}}\sum_{a_{s-1}}\frac{\phi(t_{s-1})} {p_{s-1}} \sum_{p_s}\sum_{a_s}\phi(t_s)p_1 \cdots p_{s-1}\right)
\]
\[
\le   \sum_t \frac{1}{\phi(t)}  \sumti   (2\log\log N + 2C\log t)^{s-1} \frac{N}{tt_0} 
\] 
\[
\le N  \left( \sum_t f(t)\gamma^*(t)\right) (2\log\log N + 2C\log t)^{s-1}
\]
\[
\le N \sum_{0 \le r \le s-1} \sum_t f(t)\gamma^*(t)\binom{s-1}{r}(2\log\log N)^r(2C\log t)^{s-1-r}
\]
\[
\le N \sum_{0 \le r \le s-1} (2C)^s\binom{s-1}{r} \left( \sum_t f(t) \gamma^*(t) (\log t )^s\right) (\log\log N)^{s-1}
\]
\[
\le NG(\log\log N)^{s-1}
\]
for some positive constant $G$ depending only on $s$ (cf., Lemma \ref{D})
\end{proof}

We now give our upper bound for $\osn$.

\begin{theorem}\label{upperthm}
There is a constant $F$ depending only on $s$ with 
\[
s!\osn \le \frac{N}2 (\sum_t f(t) \gamma^*(t)) (\log\log N)^s + N F (\log \log N)^{s-1}.
\]
\end{theorem}
\begin{proof}
For any $i \le s$ we have 
\[
\phi(t_i) \sum_{p_i} \sum_{a_i} \left( \frac{1}{p_i^{a_i}} - \frac{1}{p_i^{1+a_{i}}} \right) 
\]
\[
 \le \phi(t_i) \sum_{p_i} \sum_{j=1}^\infty \left( \frac{1}{p_i^{j}} - \frac{1}{p_i^{j+1}} \right)
 \le \log\log N + C \log t
\] 
by Lemma \ref{carl}.  Applying this inequality and that of the previous lemma to the inequality (\ref{5.2}) we deduce that
$s!\osn$ is at most
\[
\frac{N}2\sum_t f(t)\sumti \sum_{\langle p_i,a_i\rangle}  \sum_{p_s}\sum_{a_s} \prod_{i=1}^{s}\phi(t_i)\left(\frac{1}{p_i^{a_i}} - \frac{1}{p_i^{1+a_i}} \right)
 +   2^sGN(\log\log N)^{s-1}
\]
\[
\le \frac{N}2 \sum_t f(t)\gamma^*(t) (\log \log N +C\log t)^s +  2^sGN(\log\log N)^{s-1}
\]
\[
\le \frac{N}2(\sum_t f(t)\gamma^*(t))(\log \log N)^s  
\]
\[ 
+\frac{N}2\sum_t f(t)\gamma^*(t) \sum_{r=1}^s \binom{s}{r}(\log\log N )^{s-r} (C\log t)^r   + 2^sGN(\log\log N)^{s-1}
\]
\[
\le \frac{N}2(\sum_t f(t)\gamma^*(t))(\log \log N)^s  
\]
\[ 
+ \frac{N}2(\sum_t f(t)\gamma^*(t)(\log t)^s)(\sum_{r=1}^s \binom{s}{r}C^r)(\log\log N )^{s-1} +2^sGN(\log\log N)^{s-1}
\]
\[
\le \frac{N}2(\sum_t f(t)\gamma^*(t))(\log \log N)^s +FN(\log\log N)^{s-1}
\] 
for a constant $F>0$ depending only on $s$ (using again Lemma \ref{D}).
\end{proof}  

\section{Proof of the main theorem}\label{s: mainproof}

We will prove that for sufficiently large $N$ (depending only on the choice of $s$) the ratio
\[
\rho_s(N):= \sum_{r \le s} \mathcal{ O}_r(N) \Big/ \sum_{r \le s} \mathcal{ T}_r(N)
\]
is less than $1/2^s$.

We define $\sigma(x) = \sum f(t)\gamma^*(t)$, where the summation is over all odd $t$ with $3 \le t \le x$  (for any $x>0$). The function $\sigma(x)$ is bounded above (c.f., Corollary \ref{corD}).

The Sylow-cyclic Frobenius complements of breadth $0$ are cyclic groups and correspond to proper Frobenius triples $(m,n,\langle r + m\mathbb{Z} \rangle)$ with $m=1$,  and hence 
$\mathcal{ T}_0(N) \le N$ and $\mathcal{ O}_0(N) \le (N+1)/2$.  This, together with the theorems of the previous two sections, show that for $s>0$ we have that
\[
s!\sum_{r=0}^{s}  \mathcal{ O}_r(N) \le \frac{N}2(\log\log N)^s\sigma(\sqrt{N}) + JN(\log\log N)^{s-1},
\]
and that $s!\sum_{r=0}^{s} \mathcal{ T}_r(N)$ is at least
\[
  N\left(   \sum_{2 \le t \le \log\log N}f(t)\gamma^*(t)\right) (\log\log N)^s  \; - \; HN(\log\log N)^{s-1}\log\log\log N
\] 
for constants $J$  and $H$ depending only on $s$.  

One checks that 
\begin{equation}\label{sumest}
\sum_{r=1}^{4} \frac{r^s}{4^r} \ge \frac{2^s}{6} 
\end{equation}
for $s \le4$.  Indeed the equation (\ref{sumest}) holds for all $s$ since if $s > 4$, then
\[
\sum_{r=1}^{4} \frac{r^s}{4^r} \ge \frac{2^s}{4^2} + \frac{3^s}{4^3} \ge 2^s(\frac{1}{16}+ \frac{(3/2)^5}{64} ) \ge \frac{2^s}{6}.
\]

Thus by Lemma  \ref{Lem3.1} (C) we have
 \[
\sum_{r=0}^{4} \frac{\gamma ( 2^r)}{4^r} = \sum_{r=0}^{4}\frac{ (r+1)^s-r^s}{4^r} = \frac{5^s}{256} + 3\sum_{r=1}^{4} \frac{r^s}{4^r} \ge \frac{1}{60} + 2^{s-1}.
\]
Consequently using Lemma \ref{Lem3.1} (B)   we deduce that 
\[
\sum_{2 \le t \le \log\log N} f(t)\gamma^*(t) \ge \sum_{r=0}^{4} \quad \sum_{ \substack{ 3 \le t \le (\log\log N)/16\\2\nmid t}}f(2^r t)\gamma^*(2^rt)
\]
\[
\ge \sum_{r=0}^{4} \quad  \sum_{ \substack{ 3 \le t \le (\log\log N)/16\\2\nmid t}} \frac{1}{ (2^rt)^2\phi(t_0)}\gamma(2^r)\gamma^*(t)  
\]
\[
=\sum_{r=0}^{4}\frac{\gamma(2^r)}{4^r} \sigma((\log\log N)/16)  \ge (\frac{1}{60} + 2^{s-1} ) \sigma((\log\log N)/16).
\]

Since $\sigma(\sqrt{N}) > \sigma(3)$ therefore the ratio $\rho_s(N)$ is at most
\[
\frac{(N/2)(\log\log N)^s\sigma(\sqrt{N}) + J N (\log\log N)^{s-1}}{N(\frac{1}{60} + 2^{s-1} )\sigma((\log\log N)/16)(\log\log N)^s-HN(\log\log N)^{s-1}\log\log\log N}
\]
\[
\le \frac{1 + \frac{2J}{\sigma(3)\log\log N}}{  (\frac{1}{30} + 2^{s }) \frac{\sigma((\log\log N)/16)}{\sigma(\sqrt{N})} - \frac{2H\log\log\log N}{\sigma(3) \log\log N}}.\]

Now let $  \epsilon =  \frac{1}{60(2^s) + 31}$.  We can pick $N$ sufficiently large (where ``sufficiently large" depends only on $s$) that the numerator in the last display is at most $1+\epsilon$,  that $\frac{2H\log\log\log N}{\sigma(3) \log\log N}$ is less than $\epsilon$ and that 
$\frac{\sigma((\log\log N)/16)}{\sigma(\sqrt{N})}$ is greater than $1-\epsilon$.  
(This is possible since $\sigma$ is nondecreasing and bounded above.)  Then one checks that 
\[
\rho_s(N) < \frac{1+\epsilon}{(\frac{1}{30} + 2^{s }) (1-\epsilon) - \epsilon} = \frac{1}{2^{s}},
\]
which was to be proven.

%\begin{remark}
 %\end{remark}

\section{Upper and lower bounds for $\tsn$, $\osn$ and $\rho_s(N)$ }\label{s: morebounds}

We adapt the arguments of  Sections  \ref{s: lower} and \ref{upper} to give upper  bounds for   $\tsn$  and lower bounds for $\osn$. We will then be able to give the exact
proportion of  Sylow-cyclic Frobenius complements of width at most $s$ which have odd order  in terms of the infinite sums of  Corollary \ref{corD}.  
Rather crude  estimates of these infinite sums will show that the above proportion is between $1/(6\Upsilon)$ and $1/((21 + 6(2^s))\Upsilon)$ where
$\Upsilon = \sum_{i=1}^{\infty} \frac{i^s}{4^i}$.

Using Proposition \ref{prop2} and the arguments of Section \ref{s: lower} for a lower bound for $\tsn$, we easily obtain a lower bound for $\osn$: 
\[
s!\osn \ge \frac{N}{2} \sigma(\log\log N) (\log\log N)^s - N I' (\log\log N)^{s-1} \log\log\log N
\]
for some constant $I'$ depending only on $s$. Of course this implies that

\[
s!\sum_{r=0}^{s}  \mathcal{ O}_r(N) \ge \frac{N}{2} \sigma(\log\log N) (\log\log N)^s - N J' (\log\log N)^{s-1} \log\log\log N
\]
for a constant $J'$ depending only on $s$. Combining our upper and lower bounds on $\osn$ (and also Corollary \ref{corD}) we can deduce  that 
\begin{equation}\label{lowerosn}
s!\sum_{r=0}^{s}  \mathcal{ O}_r(N)\thicksim \frac{N}{2} \sigma( N) (\log\log N)^s
\end{equation}
(i.e., the ratio of the two sides approaches $1$).

Similarly Proposition \ref{prop1} and the arguments of 
Section \ref{upper} for an upper bound for $\osn$ easily adapt to give an upper bound for $\tsn$, 
namely,
\[
s!\tsn \le N (\log\log N)^s (\sum_{2 \le t \le \sqrt{N}}  f(t)\gamma^*(t)) + N F' (\log\log N)^{s-1}  
\]
for a constant $F'$ depending only on $s$, so that

\[
s!\sum_{r=0}^{s}  \mathcal{ T}_r(N) \le N (\log\log N)^s (\sum_{2 \le t \le \sqrt{N}}  f(t)\gamma^*(t)) + N H' (\log\log N)^{s-1}  
\]
for a constant $H'$ depending only on $s$, and hence
\begin{equation}\label{uppertsn}
s!\sum_{r=0}^{s}  \mathcal{ T}_r(N) \thicksim N (\log\log N)^s (\sum_{2 \le t \le N}f(t)\gamma^*(t)).
\end{equation}

Combining the above formulas     (\ref{lowerosn}) and (\ref{uppertsn}) with the convergence of the infinite sums in Corollary \ref{corD}, we have the

\begin{theorem}\label{exact}
\[
\lim \rho_s(N) = \frac{\sum_{ \substack { 2 \le t < \infty\\2\nmid t}} f(t)\gamma^*(t)}{2\sum_{   2 \le t < \infty } f(t)\gamma^*(t)}.
\]
\end{theorem}

Recall that  $\Upsilon = \sum_{r=1}^{\infty} r^s/4^r$.

\begin{theorem}\label{supinf}
 
\[
\frac{1}{6\Upsilon} \ge \lim  \rho_s(N)   \ge \frac{1}{ \Upsilon (21 + 6(2^s)) }.
\]
\end{theorem}

\begin{proof}

The   left-hand inequality of our theorem follows from the previous theorem and the calculation (using Lemma 3.1 (B) and (C)):

\[
\sum_{t=2}^{\infty}  f(t)\gamma^*(t) \ge \sum_{r=0}^{\infty} \sum_ {\substack{ t=3 \\ 2 \nmid t} }^{\infty} f(2^rt)\gamma^*(2^rt)
\ge \sum_{r=0}^{\infty} \quad \sum_ {\substack{ t=3 \\ 2 \nmid t} }^{\infty} \frac{1}{4^rt^2\phi(t_0)}\gamma(2^r)\gamma^*(t)
\]
\[
=   \sum_{r=0}^{\infty}  \frac{(r+1)^s-r^s}{4^r} \quad \sum_ {\substack{ t=3 \\ 2 \nmid t} }^{\infty} f(t)\gamma^*(t)
 =   3\Upsilon \sum_ {\substack{ t=3 \\ 2 \nmid t} }^{\infty} f(t)\gamma^*(t) .
\]

Each element of $\Gamma(t)$ has an entry $1$ in, say, $r$ coordinates where $0 \le r < s$, and for each such element the remaining $s-r$ coordinates can be filled with positive integers in at most $\gamma^*(t)$
ways.  (Replacing the $1$'s in the $s$-tuple by $t$'s would give an element of $\Gamma^*(t)$.)  Thus 
\begin{equation}\label{gammabdd}
\gamma(t) \le \sum_{0 \le r \le s-1} \binom{s}{r} \gamma^*(t) \le (2^s -1)\gamma^*(t).
\end{equation}

The right-hand inequality of our theorem follows from the following calculation using the above inequality (\ref{gammabdd}) and   Lemma \ref{Lem3.1} (A)
 and  (D):

%Setting $\gamma^*(1) = 1 = \phi(1)$) we deduce that  
\[
\sum_{t=2}^{\infty}  f(t)\gamma^*(t) \le \sum_{r=0}^{\infty} \sum_ {\substack{ t=3 \\ 2 \nmid t} }^{\infty} f(2^rt)\gamma^*(2^rt)
+ \sum_{r=1}^{\infty} f(2^r)\gamma^*(2^r)
\]
\[
\le \sum_{r=0}^{\infty} \sum_ {\substack{t=3 \\ 2 \nmid t} }^{\infty} \frac{1}{4^rt^2\phi(t_0)}\gamma(2^r)\gamma(t) + \sum_{r=1}^{\infty}\frac{r^s-(r-1)^s}{4^r}
\]
\[
\le 3\Upsilon (2^s -1)  \sum_{ \substack { t=3  \\2\nmid t}}^{\infty} f(t)\gamma^*(t) + \frac34\Upsilon
\]
\[
 \le (3(2^s -1) + \frac{3/4}{1/18}) \Upsilon  \sum_{ \substack { t=3\\2\nmid t}}^{\infty} f(t)\gamma^*(t)  ,
\]
since $\sum_{ \substack { t=3\\2\nmid t}}^{\infty} f(t)\gamma^*(t) > f(3)\gamma^*(3) = 1/18$.

\end{proof}

\begin{remark} For small $s$ the left-hand inequality of the last theorem gives a modest improvement of the Main Theorem, but  the improvement can be substantial for large $s$.  For example,  when  $s=10$ then $6\Upsilon$ is about $99,851$, almost $100$ times the size of $2^{10}$. 
Indeed
if $s=20$ then $6\Upsilon$ exceeds $25 \times 10^{14}$ while $2^s$ is of course a bit more than a million.  (The calculations of $\Upsilon$ were made using WolframAlpha.)
\end{remark}

\bibliographystyle{plain}
\bibliography{\jobname} 

\begin{thebibliography}{00} 

\bibitem{A}
 M. Aschbacher, 	\textit { Finite Group Theory, 2nd ed.}. (Cambridge University Press, Cambridge, 2000).
  \bibitem{Br}
 R. Brown, Frobenius groups and classical maximal orders,  	\textit { Mem. Amer. Math. Soc.} {\bf 717} (2001).

\bibitem{DM}
 J. Dixon and B. Mortimer, 	\textit { Permutation Groups}. (Springer-Verlag, New York, 1996).
 

 \bibitem{HW}
G. H. Hardy and E. M. Wright,  	\textit {An Introduction to the Theory of Numbers, 5th ed.}. (Oxford University Press, Oxford, 1979).

 \bibitem{L}
 T. Y. Lam, Finite groups embeddable in division rings, 	\textit { Proc. Amer. Math. Soc.}  {\bf 129} (2001), 3161-3166.

 \bibitem{BT}
 H. L. Montgomery and R. C. Vaughn, The large sieve,	\textit { Mathematika} {\bf 20} (1973), 119-134.

% \bibitem{P}
 %C. Pomerance, On the distribution of amicable numbers,	\textit {J. Reine Agnew. Math.} {\bf 293/294} (1977), 217-222.

 \bibitem{N}
K. K. Norton, On the number of restricted prime factors of an integer,	\textit {Illinois J. Math.} {\bf 20} (1976), 681-705.



\end{thebibliography}

\end{document}